\documentclass{amsart}
\usepackage{hyperref}
\usepackage{amssymb}
\usepackage{fancyhdr}		
\usepackage{color}

\usepackage[T1]{fontenc}
\usepackage[latin1]{inputenc}
\usepackage{url}	
\usepackage{graphicx}	
\usepackage{amssymb}
\usepackage{amsmath}
\usepackage{amsthm}
\usepackage[T1]{fontenc}
\usepackage{ae}

\newtheorem{theo}{Theorem}[section]
\newtheorem{cor}[theo]{Corollary }
\newtheorem{lemm}[theo]{Lemma}

\newtheorem{defi}{Definition}[section]

\newtheorem{prop}[theo]{Proposition}

\title[]{Asymptotic behavior of Moncrief Lines in constant curvature space-times}

\date{\today}
\author{M. Belraouti}
\address{Mehdi Belraouti \newline
Facult\'e de Math\'ematiques,\\
USTHB, BP 32, El-Alia,\\
16111 Bab-Ezzouar, Alger (Algeria)}
\email{mbelraouti@usthb.dz}

\author{A. Mesbah}
\address{Abderrahim Mesbah \newline
Department of Mathematics,\\
University of Luxembourg,\\
}
\email{abderrahim.mesbah@uni.lu}

\author{L. Messaci} 

\address{Lamine Messaci\newline
Laboratoire Jean Alexandre Dieudonn\'e,\\
Parc Valrose, 28 Avenue Valrose\\
06108-Nice, France}
\email{messaci@unice.fr}

\begin{document}
\maketitle

\noindent{\bf Abstract.}
We study the asymptotic behavior of Moncrief lines on $2+1$ maximal globally hyperbolic spatially compact space-time $M$ of non-negative constant curvature. We show that when the unique geodesic lamination associated with $M$ is either maximal uniquely ergodic or simplicial, the Moncrief line converges, as time goes to zero, to a unique point in the Thurston boundary of the Teichm\"uller space.

\tableofcontents

\section{Introduction}

A space-time is an oriented chronologically oriented Lorentzian manifold. It is said to be globally hyperbolic (we write $GH$) if it admits a Cauchy hypersurface, which is a spacelike hypersurface  such that every inextensible causal curve (curve for which the norms of the tangent vectors are non-positive) intersects it at exactly one point. If in addition the Cauchy hypersurface is compact, then we say that the space-time is globally hyperbolic spatially compact and we write $GHC$. Another equivalent way to define global hyperbolicity is through the existence of a proper function, called a Cauchy time function, which is strictly increasing along causal curves and surjective on the inextensible causal ones. The level sets of such a function are Cauchy hypersurfaces, and the spacetime is globally hyperbolic spatially compact. By a classical result of Geroch \cite{Geroch30}, every $GH$ space-time is diffeomorphic to the product of a Cauchy  hypersurface $S$ by $\mathbb{R}$.

In general relativity, the Einstein equation governs how the distribution of matter is related to the metric properties of the space-time. However, it does not, on its own, determine the space-time and the main problem of gravity is to classify its solutions. A classical way to deal with this problem is to solve the associated Cauchy problem. Namely, this consists in considering an initial data $(S, h, \amalg)$ of a Riemannian metric together with a symmetric $2$-tensor defined on a manifold $S$. Then trying to find a Lorentzian metric $g$ on $M=S\times\mathbb{R}$ verifying the Einstein equation and such that $h$ and $\amalg$ are respectively the Riemannian metric and the second fundamental form defined on $S\times\left\{0\right\}$ by $g$. In dimension $2+1$, space-times of constant curvature are the unique solutions of the Einstein equation and hence for the Cauchy problem to admit a solution, the initial data $(S, h, \amalg)$ must satisfy the classical Gauss-Codazzi equations. Any initial data satisfying the Gauss-Codazzi equations gives rise to a  solution of the Einstein equation \cite{Geroch9}. Moreover, the solution is globally hyperbolic $GH$.  
By Choquet-Bruhat-Geroch \cite{Geroch9} this local solution has a unique maximal  extension. In this sense, a globally hyperbolic space-time $(M,g)$ of constant curvature is said to be maximal if it does not extend to a constant curvature globally hyperbolic space-time which is also solution of the Einstein equation. The  Cauchy problem is then reduced to the classification of $MGHC$ space-times of constant curvature.

In \cite{mess1}, using geometric methods, Mess gives a full classification of $2+1$ maximal globally hyperbolic spatially compact ($MGHC$) space-times of non-positive constant curvature. This classification was completed in the works of Scannell, Barbot, Beguin, Zeghib, Bonsante and Benedetti for all constant curvature and all dimension cases (\cite{scannell1}, \cite{barbot1}, \cite{bonsante1}, \cite{BBZ1}). In particular, there is a one to one correspondance between measured geodesic laminations on a closed surface $S$ and MGHC space-times of dimension $2+1$ admitting a Cauchy surface diffeomorphic to $S$.

A different approch consists in reducing the dimension of the space of solutions by imposing some additional contraints i.e choosing a gauge. 
In dimension $2+1$, one way to do that is the $CMC$ reduction in which every spacelike surface $S\times\left\{*\right\}$ is supposed to have constant mean curvature. Thanks to the works of Andersson, Barbot, Beguin and Zeghib \cite{bbz4}, \cite{bbz5}, \cite{BBZ1}, every $2+1$ $MGHC$ space-time of constant curvature posses a time function, called $CMC$ time function, where the levels have constant mean curvature. Therefore, 
the $CMC$ reduction  applies in these cases.

Using $CMC$ reduction, Moncrief (\cite{Moncrief2}) proved that the Einstein equation reduce to a non-autonomous Hamiltonian system, called the Moncrief flow, on the cotangent bundle of the Teichm\"uller space of $S$. In particular, to every orbit 
$\alpha$ of this flow corresponds a unique constant curvature $MGHC$ space-time of dimension $2+1$ such that the projection curve of $\alpha$ on $\operatorname{Teich}(S)$, that we call the Moncrief line,  is the curve $\left([g_{a}^{T_{cmc}}]\right)_{a}$ of conformal classes  of the Riemannian metrics $g_{a}^{T_{cmc}}$  defined on $S$ by the $CMC$ time $T_{cmc}$. This allows us  to study their asymptotic behavior with respect to the time parameter. It is a natural question that one can formulate for general time function: Let $T$ be a Cauchy time function defined on a non-elementary $2+1$ maximal globally hyperbolic spatially compact space-time $M=S\times \mathbb{R}$ of constant curvature. Such a function defines naturally a one parameter family of Riemannian metrics on $S$ which in turn gives rise (by considering their conformal classes) to a curve $\left([g^{T}_ {a}]\right)_{a}$, parameterized by time $T$, in the Teichm\"uller space $\operatorname{Teich}(S)$. 

Benedetti and Bonsante (\cite{benebonsan1}), proved that the curve corresponding to the cosmological time $T_{cos}$ is nothing but the grafting line $(\operatorname{gra}_{\frac{\lambda}{a}}(S))_{a}$ defined by the measured geodesic lamination $\lambda$ associated to $M$. It is real analytic and converges, as
time  goes to $+\infty$, to the hyperbolic
structure $S$. In the flat case, Belraouti (see \cite{mehdi3}) proved  that  near infinity, the Moncrief lines and the Teichm\"uller curves defined by the $k$-times (i.e time functions where the levels have constant Gauss curvature) behave in the same way as the grafting lines.

Near $0$ this is no larger the case when the $k$-time is concerned. Indeed, on the one hand, It is proved in \cite{mehdi3,mehdi1} that the curve $([g^{T_{k}}_{a}])_{a>0}$  converges when time goes to $0$ to the point, in the Thurston boundary, corresponding to the measured geodesic lamination $\lambda$. On the other hand, Diaz and Kim \cite{Kim} showed that the grafting line defined by a simplicial measured geodesic lamination $\sum c_{i}\gamma_{i}$ converges to the point, in the Thurston boundary of the Teichm\"uller space, corresponding to the measured geodesic lamination $\sum \gamma_{i}$.

In this article we will study the asymptotic behavior of the Moncrief lines when time goes to $0$. We will show that for a maximal uniquely ergodic or a simplicial geodesic lamination $\lambda$, the Moncrief lines and the grafting lines have the same limite when time goes to $0$. More precisely:

\begin{theo}
\label{theoprinc}
Let $S$ be a closed hyperbolic surface and let $\lambda$ be a measured geodesic lamination on $S$. Let $M$ be the unique $2+1$ non elementary $MGHC$  space-time of non negative constant curvature associated to $\lambda$. Let $T_{cmc}$ be  the $CMC$ time of $M$.
\begin{enumerate}
\item If $\lambda$ is maximal uniquely ergodic then the curve $([g^{T_{cmc}}_{a}])_{a>0}$ converges to the projective lamination $\left[\lambda\right]$ in  the Thurston boundary. 
\item If $\lambda=\sum c_{i}\gamma_{i}$ is simplicial then the curve  $([g^{T_{cmc}}_{a}])_{a>0}$ converges to the projective lamination $\left[\sum \gamma_{i}\right]$ in  the Thurston boundary.
\end{enumerate}
\end{theo}

\section{Generalities}
\subsection{Teichm\"uller distance, length spectrum and Thurston's compactification}
Let  $S$ be a closed surface of genus $g \geq 2$. Denote by $\operatorname{Teich}(S)$ the Teichm\"uller space of $S$. The Teichm\"uller distance between two points $h,h' \in \operatorname{Teich}(S)$ is defined by $$d_{\operatorname{Teich}(S)}(h,h')=\frac{1}{2}\operatorname{min}_{f}\log K(f),$$ where $K(f)$ is the optimal quasi-conformal constant of $f$, and the minimum is taken over all quasi-conformal homemorphisms $f:(S,h) \to (S,h')$ isotopic to the identity. The quasi-conformal map that realizes this minimum is unique and is known as the extremal map.

Let $\mathcal{S}$ be the set of free homotopy classes of simple closed curves non-contractible on $S$. For every $\alpha\in \mathcal{S}$, and every $h\in \operatorname{Teich}(S)$ we denote by $l_{h}(\alpha)$ the hyperbolic length of the unique simple closed geodesic in the free homotopy class $\alpha$ with respect to $h$. The length of a free homotopy class $\alpha\in \mathcal{S}$ with respect to two elements $h,h' \in \operatorname{Teich}(S)$ is related to their Teichm\"uller distance. Namely, we have:
\begin{prop} {\em (\cite[Lemma~3.1]{Wolpert})}\label{proppropprop}
 $\forall \alpha \in \mathcal{S}$:
 $$e^{-2d_{\operatorname{Teich}(S)}(h,h')}\ell_{h}(\alpha) \leq \ell_{h'}(\alpha) \leq e^{2d_{\operatorname{Teich}(S)}(h,h')}\ell_{h}(\alpha)$$
\end{prop}

For every $\alpha \in \mathcal{S}$, the map $h \mapsto l_h(\alpha)$ is continuous (see \cite{fathi2021thurston}). This gives rise to the continuous length functional
$$
l: \operatorname{Teich}(S) \longrightarrow \mathbb{R}_{+}^{\mathcal{S}}
$$
defined by
$$
l: h \mapsto l_h: \alpha \mapsto l_h(\alpha),
$$
which, when composed with the projection
$$
\pi: \mathbb{R}_{+}^{\mathcal{S}} \longrightarrow \mathbb{P}(\mathbb{R}_{+}^{\mathcal{S}}),
$$
turns out to be an embedding whose image is relatively compact (see \cite{fathi2021thurston}). Moreover, Thurston (see \cite{fathi2021thurston}) showed that the complement of this image in its closure coincides with the image of the projectivization $\mathcal{P}\mathcal{M}\mathcal{L}(S)$ of the space of measured geodesic laminations $\mathcal{M}\mathcal{L}(S)$ on $S$, via the embedding defined by intersection number functional $i$. Additionally, $\overline{\operatorname{Teich}(S)}=\operatorname{Teich}(S)\sqcup \mathcal{P}\mathcal{M}\mathcal{L}(S)$ is homeomorphic to a closed $6g-6$-dimensional ball. This compactification is known as the Thurston compactification of the Teichm\"uller space $\operatorname{Teich}(S)$. Thus, a sequence $(h_{n})_{n}$ in $\operatorname{Teich}(S)$ converges to the projective class $[\lambda]\in \mathcal{P}\mathcal{M}\mathcal{L}(S)$ if there is a sequence $(\epsilon_{n})_{n}$ of positive real numbers such that for all $\alpha\in \mathcal{S}$, we have $\lim_{n \to \infty}\epsilon_{n}l_{h_{n}}(\alpha)=i(\lambda,\alpha)$.

\subsection{Grafting}
For details see \cite{kulkarni1994canonical}\cite{kamishima1992deformation} \cite{dumas2009complex}. Let $\mathcal{P}(S)$ denote the space of complex projective structures on $S$. Recall that a complex projective structure on $S$ is a maximal atlas where the charts are modeled on $\mathbb{CP}^1$ and have M\"obius transition functions. Using the grafting, Thurston provides a geometric description of the space of complex projective structures  by showing that $\mathcal{P}(S)$ is homeomorphic, via the grafting map $ \operatorname{Gr}: \operatorname{Teich}(S)\times \mathcal{ML}(S) \to \mathcal{P}(S)$, to $\operatorname{Teich}(S)\times \mathcal{M}\mathcal{L}(S)$. The grafting map $\operatorname{Gr}: \operatorname{Teich}(S) \times \mathcal{ML}(S) \to \mathcal{P}(S)$ can be described as follow: let $h$ be a hyperbolic metric on $S$, and let $\gamma$ be a simple closed geodesic on $S$ with a weight $t$. Define $X^{\gamma}_{t} := \gamma \times \left [ 0,t \right ]$ to be an Euclidean cylinder and replace $\gamma$ with $X^{\gamma}_{t}$. Then we get a surface with a $C^{1,1}$ metric (Thurston's metric) which is hyperbolic on $S \setminus \gamma$ and euclidean on $X^{\gamma}_{t}$. Then $\operatorname{Gr}_{\gamma}h$ is the canonical projective structure  piecing together the original projective structures of $h$ and $X^{\gamma}_{t}$. Grafting extends by continuity 
to general measured laminations. The conformal grafting is the map $ \operatorname{gr}: \operatorname{Teich}(S)\times \mathcal{ML}(S) \to \operatorname{Teich}(S)$ which associates to a pair $(h,\gamma)$ the underlying conformal structure of $\operatorname{Gr}_{\gamma}h$.

\subsection{Twisting numbers}
For details see \cite{minsky1996extremal}. Let $(S,h)$ be a hyperbolic surface. The twisting number of a simple closed curve $\gamma$ along an oriented simple closed curve $\beta$ is defined as follows: let $\gamma_{h}$ and $\beta_{h}$ be the unique geodesics homotopic to $\gamma$ and $\beta$ respectively with respect to the metric $h$. For each $x \in \gamma_{h} \cap \beta_{h}$, let $\Tilde{x}$ to be a lift of $x$, and let $\Tilde{\gamma}$ and $\Tilde{\beta}$ be the respective  lifts of $\gamma$ and $\beta$ that intersect at $\Tilde{x}$. Note that the geodesic $\Tilde{\beta}$ is oriented since $\beta$ is oriented. 
Let $a_{r}$ be the endpoint of $\Tilde{\gamma}$ in $\partial_{\infty}\mathbb{H}^2$ to the right of $\Tilde{\beta}$, and let $a_{l}$ be the endpoint of $\Tilde{\gamma}$ in $\partial_{\infty}\mathbb{H}^2$ to the left of $\Tilde{\beta}$.\\
Denote respectively by $pr(a_{r})$ and $pr(a_{l})$ the projections of $a_{r}$ and $a_{l}$  on $\Tilde{\beta}$. Note that $pr(a_r) - pr(a_l)$ depends only on the intersection point $x$. Then the twisting number $Tw$ is defined to be:
$$Tw_{h}(\gamma,\beta) := \left| \min\limits_{\gamma_{h}\cap \beta_{h}} \frac{pr(a_r) - pr(a_l)}{\ell_{h}(\beta)} \right| $$

In the rest of the paper we will need  the following lemma:
\begin{lemm}
\label{bounded}
Let $h,h' \in \operatorname{Teich}(S)$, let $\gamma$ be a non-trivial simple closed curve in $S$ and let $\beta$ be an oriented simple non-trivial closed curve. For any $M > 0$ and any $K > 0$, there exists $M'>0$ depending  only on $M$ and $K$, such that if $d_{\operatorname{Teich}(S)}(h,h') < K$ and $\ell_{h}(\beta)Tw_{h}(\gamma,\beta) < M$, then $\ell_{h'}(\beta)Tw_{h'}(\gamma,\beta) < M'$. 
\end{lemm}
\begin{proof}

Let $\gamma_{h}$ and $\beta_{h}$ be the unique simple closed geodesics homotopic to $\gamma$ and $\beta$  with respect to the metric $h$. Similarly, let $\gamma_{h'}$ and $\beta_{h'}$ be the unique simple closed geodesics homotopic to $\gamma$ and $\beta$  with respect to the metric $h'$.

Let $f:(S,h) \to (S,h')$ be the extremal map and let $\Tilde{f}: \mathbb{H}^2 \to \mathbb{H}^2 $  be a lift of $f$, where $\mathbb{H}^2 $ is the half plane model of the hyperbolic space. Recall that $\Tilde{f}$ extends continuously to a $k$ quasi-symmetric map $\partial \Tilde{f}: \partial_{\infty} \mathbb{H}^2 \to \partial_{\infty} \mathbb{H}^2$, where $k$ depends only on $K$.

Since $\gamma_{h}$ and $\beta_{h}$ intersect at a finite number of points, there are a point $x\in \gamma_{h}\cap \beta_{h}$ along with two lifts $\Tilde{\gamma}_{h}$ and $\Tilde{\beta}_{h}$ of $\gamma_{h}$ and $\beta_{h}$, respectively, such that  $Tw_{h}(\gamma,\beta) = \left| \frac{pr(a_{r,h}) - pr(a_{l,h})}{\ell_{h}(\beta)} \right|$, where $pr(a_{r,h}), pr(a_{l,h})$ are the projections on $\Tilde{\beta}_{h}$ of the end points $a_{r,h}$, $a_{l,h}$  of $\Tilde{\gamma}_{h}$ lying respectively on the right and left side of  $\Tilde{\beta}_{h}$. 

Without loss of generality, up  to conjugacy by elements of $PSL(2,\mathbb{R})$,  assume that $\Tilde{\beta}_{h}$ is the geodesic $(0,\infty)$ oriented from $0$ to $\infty$ and one endpoint of  $\Tilde{\gamma}_{h}$,  which is on right of $\Tilde{\beta}_{h}$, is  $1$. Assume further that  
the quasi-symmetric map $\partial \Tilde{f}$ fixes the three points $ \left\{0,1,\infty \right\}$. Since $f$ is isotopic to the identity, the oriented geodesic $(0,\infty)$ must  also be a lift of $\beta_{h'}$. Moreover, the geodesic $\gamma_{h'}$ admits a lift $\Tilde{\gamma}_{h'}$ such that the end points $a_{r,h'}$, $a_{l,h'}$ of $\Tilde{\gamma}_{h'}$ to the right and left side of  $\Tilde{\beta}_{h'}$ are  given by $\partial \Tilde{f}(1)=1$ and $\partial \Tilde{f}(a_{l,h})$, respectively. Thus $\ell_{h}(\beta)Tw_{h}(\gamma,\beta) = \left| \ln(-a_{l,h}) \right|$, and $\ell_{h'}(\beta)Tw_{h'}(\gamma,\beta) \leq \left| \ln(-\partial \Tilde{f}(a_{l,h})) \right|$. 

Now, on the one hand $\ell_{h}(\beta)Tw_{h}(\gamma,\beta) < M$ implies that $-e^{M}<a_{l,h}<-e^{-M}$. On the other hand, since $\partial \Tilde{f}$ is $k-$quasi-symmetric,  there exists a constant $M'$ depending only on $k$ and $M$ such that $\partial \Tilde{f}(]-e^{M},-e^{-M}[)\subset ]-e^{M'},-e^{-M'}[$. Thus  $\ell_{h'}(\beta)Tw_{h'}(\gamma,\beta) \leq \left| \ln(-\partial \Tilde{f}(a_{l,h})) \right|< M'$.

\end{proof}
\subsection{Constant curvature space-times}
In this section we will review certain properties of constant curvature space-times.
Recall that the models of constant curvature space-times are:
\begin{enumerate}
\item The Minkowski space $\mathbb{R}^{1,n}$. That is the vectoriel space $\mathbb{R}^{n+1}$ endowed with the standard Lorentzian metric $q_{1,n}=-dx_{0}^{2}+dx_{1}^{2}+...+dx_{n}^{2}$.  It is a globally hyperbolic spactime whose isometry group is the Poincar\'e group $O(1,n)\ltimes\mathbb{R}^{1,n}$; 
\item The de Sitter space $dS_{n}$. That is the one sheeted hyperboloid $q_{1,n}=+1$ endowed with the Lorentzian metric induced by $q_{1,n}$. It is the positive curvature model space and it is a globally hyperbolic space-time whose isometry group is $O(1,n)$;
\item The anti de Sitter space $AdS_{n}$. That is the quadric $q_{2,n-1}=-1$ endowed with the Lorentzian metric induced by $q_{2,n-1}=-1$, where $q_{2,n-1}=-dx_{0}^{2}-dx_{1}^{2}+...+dx_{n}^{2}$. It is the negative curvature model space. Unlike the Minkowski and the de Sitter space-times, the anti de Sitter space-time is not globally hyperbolic. Its group of isometry is $O(2,n)$.
\end{enumerate}

It is well known that if  $(M, g)$ is a space-time with constant curvature, then every point in $M$ admits a neighbourhood which is isometric to an open subset of the corresponding constant curvature model space. This allows the Lorentzian metric on $M$ to be viewed as a $(G,X)$-structure on $M$, where in the flat case, $(G,X)=(O(1,n)\ltimes\mathbb{R}^{1,n},\mathbb{R}^{1,n})$, in the de Sitter case, $(G,X)=(O(1,n), dS_{n})$, and in the anti-de Sitter case, $(G,X)=(O(2,n), AdS_{n})$. This yields, in particular, a locally isometric developing map $\operatorname{Dev}:\tilde{M}\longrightarrow X$, as well as a holonomy representation $\rho:\pi_{1}(M)\longrightarrow G$.

Mess \cite{mess1} classified $MGHC$ space-times of dimension  $2+1$ in the flat and the anti de Sitter cases. Specifically, 
\begin{theo}
Let $(M,g)$ be a flat or anti de Sitter $MGHC$ space-time of dimension $2+1$. Then the developing map $\operatorname{Dev}:\tilde{M}\longrightarrow X$ is injective so that $M$ is isomorphic to the quotient of  the open convex domain $\operatorname{Dev}(\tilde{M})$ of $X$ by the free torsion sub-group $\rho(\pi_{1}(M))$ of $\operatorname{Isom}(X)$. Moreover, if $(M,g)$ is a non elementary flat space-time then the composition of $\rho$ with the projection $L:O(1,n)\ltimes\mathbb{R}^{1,n}\longrightarrow O(1,n)$ is injective.
\end{theo}
Following Mess work's, Barbot \cite{barbot1} extended this result to $MGH$ Cauchy complete flat space-times of any dimension.

Actually the  open convex domain $\operatorname{Dev}(\tilde{M})$ in the Mess theorem is a regular domaine. This notion was initially introduced by Bonsante in \cite{bonsante1} for Minkowski space. 
Let $\mathfrak{P}$ be the space of all lightlike hyperplanes in $\mathbb{R}^{1,n}$, $\Lambda$ be a closed subset of $\mathfrak{P}$ and consider $\Omega:=\bigcap_{P\in\Lambda} I^{+}(P)$. By \cite{barbot1}, the subset $\Omega$ is an open convex domain of $\mathbb{R}^{1,n}$. It is non empty whenever $\Lambda$ is compact. If $\Lambda$ contains more than two elements, then the open convex domain $\Omega$, if not empty, is called a future complete regular domain. A past complete regular domain is similarly defined.

Let $\Gamma$ be a torsion free uniform lattice of $SO^{+}(1,n)$. A cocycle of $\Gamma$ is an application $\tau:\Gamma\rightarrow \mathbb{R}^{1,n}$ such that $\tau(\gamma_{1}.\gamma_{2})=\gamma_{1}\tau(\gamma_{2})+\tau(\gamma_{1})$. An affine  deformation of $\Gamma$ associated to $\tau$ is the morphism $\rho_{\tau}:\Gamma\rightarrow SO^{+}(1,n)\ltimes\mathbb{R}^{1,n}$ defined by  $\rho_{\tau}(\gamma).x=\gamma.x+\tau(\gamma)$ for every $\gamma\in \Gamma$ and $x\in\mathbb{R}^{1,n}$. Bonsante \cite{bonsante1} showed that to every affine deformation of $\Gamma$ corresponds a unique (up to reorientation) future complete maximal flat regular domain $\Omega$ on which $\Gamma_{\tau}=\rho_{\tau}(\Gamma)$ acts freely properly discontinuously. Conversely, by \cite{barbot1} every future complete $MGHC$ flat non elementary space-time $M$ is up to finite cover the quotient of a future  complete regular domain by a discrete subgroup of $SO^{+}(1,n)\ltimes\mathbb{R}^{1,n}$.

Subsequently, Barbot, B\'eguin and Zeghib \cite{BBZ1} extended Mess's result to higher dimensional anti de Sitter space-times. They introduced the notion of a regular anti-de Sitter domain and demonstrated that every $n$-dimensional $MGHC$ anti-de Sitter space-time is the quotient of a regular domain by a torsion-free discrete subgroup of $O(2,n)$. 

In the de Sitter case, the developing map is no longer injective. Scannel \cite{scannell1} proved that every $MGHC$ de Sitter space-time is the quotient of a standard $dS$ space-time by a free torsion discret subgroup of $SO^{+}(1,n+1)$. More precisely, Let $d:S\rightarrow \mathbb{S}^{n}$ be a developing map of  a simply connected M\"obius manifold $S$, that is a manifold equipped with a $(G,X)$-structure, where $G=O^{+}(1,n)$ and $X=\mathbb{S}^{n}$ is the Riemannian sphere.  A round ball of $S$ is an open convex set $U$ of $S$ on which $d$ is an homeomorphism. It is said to be proper if $d(\bar{U})$ is a closed round ball of $\mathbb{S}^{n}$. Let $B(S)$ be the space of proper round ball of $S$. By a result of \cite{BBZ1}, there is a natural topology on $B(S)$ making it locally homeomorphic to $\mathbb{D}\mathbb{S}_{n+1}$. By \cite{BBZ1}, the space $B(S)$ endowed with the pull back metric of $\mathbb{D}\mathbb{S}_{n+1}$ is a simply connected future complete globally hyperbolic locally de Sitter space-time called $dS$-standard space-time.

In the $2+1$ dimensional case, Mess \cite{mess1} provides a parametrization of $MGHC$ flat non elementary space-times in terms of measured geodesic laminations. More precisely, if $S$ is a compact surface of hyperbolic type, then to every measured geodesic lamination on $S$ corresponds a unique flat $MGHC$ Lorentz structure on $S\times \left]0,+\infty \right[$. Conversely, every $MGHC$ flat space-time $M$ admitting a Cauchy surface diffeomorphic to $S$, gives rise to unique measured geodesic lamination on $S$. Using Wick rotations, Benedetti and Bonsante \cite{benebonsan1} extend this correpondence to de Sitter and anti de Sitter space-times. From now until the end of the paper we will restrict ourselves to the dimension $2+1$.

\section{Geometric time functions}
One of the specificities of constant curvature space-times is the possession of geometrically interesting time functions, namely the cosmological time, the $CMC$ time, and the $k$-time. These, in particular, provide important examples of quasi-concave time functions. 

Let $S$ be a $C^{2}$ Cauchy surface of a $MGHC$ space-time $M$ of constant curvature and consider   $\Pi_{S}$ to be its second fundamental form defined by $\Pi_{S}(X,Y)=\langle\triangledown_{X}\mathfrak{n},Y\rangle$, where $\mathfrak{n} $ is the future oriented normal vector field. The mean curvature $H_{S}$ at a point $p$ of $S$ is defined by $H_{S}=\frac{tr(\Pi)}{n}=\frac{\lambda_{1}+\lambda_{2}}{2}$, where  $\lambda_{1}, \lambda_{2}$ are the principal curvatures of $S$. The \textbf{intrinsic} curvature $k_{S}$ at a point $p$ of $S$ is defined to be $k_{S}=-det(\Pi)=-\lambda_{1}\lambda_{2}$.\\
The surface $S$ is said to be convex if its second fundamental form is positive-definite, or equivalently, if its principal curvatures are negative. The convexity of $S$ is equivalent to the geodesic convexity of $J^{+}(S)$. This enables us to extend the notion of convexity to non-smooth hypersurfaces.

\begin{defi}
A  Cauchy time function $T:M \rightarrow \mathbb{R}$ is quasi-concave if its levels are convex. We will denote by $S^{T}_{a}$ the level sets of $T$.
\end{defi}

\subsection{Cosmological time}
The cosmological time $T_{cos}$ is defined at a point $p$ by: $$T_{cos}(p)=\sup_{\alpha}\int\sqrt{-\left|\dot{\alpha}(s)\right|^{2}}$$ where the supremum is taken over all the past causal curves starting at $p$.

It gives  a simple and important first example of time functions. By \cite{bonsante1}, \cite{benebonsan1} the cosmological time of a $MGHC$ spcae-time of constant curvature is a $C^{1}$ Cauchy time function.

In the flat case, the cosmological time is concave and hence quasi-concave Cauchy time function (see \cite{bonsante1}). On the other hand, By \cite{scannell1}, \cite{BBZ1} the cosmological time is a quasi-concave time in the de Sitter case. In contrast, this is no longer true in the anti-de Sitter case.

\subsection{k-time}
The $k$-time is a Cauchy time $T_{k}:M\rightarrow \mathbb{R}$ such that every level $S^{T}_{a}$, if not empty, is of constant Gauss curvature $a$. Barbot, B\'eguin and Zeghib \cite{BBZ2} proved the existence and uniqueness of such time in the flat and de Sitter case. More precisely,
\begin{theo}
Let $M$ be a $2+1$ non elementary $MGHC$ space-time of constant curvature. Assume in addition that $M$ is future complete in both the flat and de Sitter cases. Then:
\begin{enumerate}
\item In the flat case: $M$ admits a unique $k$-time $T_{k}:M\rightarrow \left]-\infty, 0 \right[ $ ranging over $\left]-\infty, 0 \right[$;
\item In the de Sitter case: $M$ admits a unique $k$-time $T_{k}: M\rightarrow \left]-\infty, -1 \right[ $ ranging over $\left]-\infty, -1 \right[$;
\item In the anti de Sitter case: there is no  globally defined $k$-time. However, the two connected components \textbf{of the complement of the convex core} admit a unique $k$-time \textbf{ranging over} $\left]-\infty, 0 \right[$.
\end{enumerate}
\end{theo}

By definition, the $k$-time is quasi-concave.

\subsection{CMC time}
The $CMC$ time is a Cauchy time $T_{cmc}: M\rightarrow \mathbb{R}$ such that every level $S^{T}_{a}$, if not empty, is of constant mean curvature $a$. 

Andersson, Barbot, B\'eguin and Zeghib studied the existence and uniqueness of the $CMC$ time in \cite{anderson2}, \cite{anderson6}, \cite{bbz4}, \cite{bbz5}, \cite{BBZ1}. They particularly demonstrated that a $2+1$ non-elementary $MGHC$ space-time of constant curvature admits a unique $CMC$ time, which ranges over different intervals depending on the curvature type: $\left]-\infty, 0 \right[$ in the flat case, $\left]-\infty, -1 \right[$ in the de Sitter case, and $\left]-\infty, +\infty \right[$ in the anti-de Sitter case. 

In the non-negative constant curvature case, one of the interesting properties of $CMC$ time is its comparability to cosmological time. More precisely:
\begin{prop}
\label{Propoimp}
Let $M$ be a $2+1$ non elementary future complete $MGHC$ space-time of non negative constant curvature. Let $T_{cos}$, $T_{cmc}$ be respectively the cosmological time and the $CMC$ time of $M$. Then, for every $a\in ]-\infty,-1[$:
$$\frac{\sup_{T_{cmc}^{-1}(a)}T_{cos}}{\inf_{T_{cmc}^{-1}(a)}T_{cos}}\leq 2 .$$
\end{prop}
\begin{proof}
$1)$ The flat case. This is a direct consequence of \cite[Theorem~1.4]{BBZ1}.\\
$2)$ The de Sitter case. By \cite[Theorem~16.1]{BBZ1}, for every  $a\in ]-\infty,-1[$, we have:
$$H_{S^{T_{cos}}_{\frac{1}{2}\text{arcoth}(-a)}}\leq a\leq H_{S^{T_{cos}}_{\text{arcoth}(-a)}}.$$
This means that the couple of  cosmological levels $\left(S^{T_{cos}}_{\text{arcoth}(-a)} , S^{T_{cos}}_{\frac{1}{2}\text{arcoth}(-a)}\right) $ forms an $a-$barrier in the sens of \cite[Definition~4.1]{BBZ1}. Therefore, by \cite[Theorem~4.9]{BBZ1}, the $CMC$ level $S^{T_{cmc}}_{a}$ lies in the future of $S^{T_{cos}}_{\frac{1}{2}\text{arcoth}(-a)}$ and in the past of $S^{T_{cos}}_{\text{arcoth}(-a)}$. The result follows immediately.
\end{proof}
By  composing with the reparametrizations $b\longmapsto -\frac{1}{b}$ in the flat case and  $b\longmapsto \text{arcoth}(-b)$ in the de Sitter case, we assume that the $CMC$ time takes its values in $\mathbb{R}_{+}^{*}$. In this case, Proposition \ref{Propoimp} becomes:
\begin{prop}
Let $M$ be a $2+1$ non elementary future complete $MGHC$ space-time of non negative constant curvature. Let $T_{cos}$, $T_{cmc}$ be respectively the cosmological time and the $CMC$ time of $M$. Then, for every $a\in ]0,\infty[$:\\
$$T_{cos} \leq T_{cmc} \leq 2T_{cos}.$$

\end{prop}

By a result of Treibergs \cite{triberg}, the $CMC$ time is quasi-concave in the flat case. Our next Proposition shows that this is also true in the de Sitter case at least in the $2+1$ dimensional case: 
\begin{prop}
Let $M$ be a $2+1$ non elementary future complete $MGHC$ space-time of positive constant curvature. Then the $CMC$ time $T_{cmc}$ of $M$ is quasi-concave.
\end{prop}
\begin{proof}
Let us fix $a\in ]-\infty,-1[$. We need to prove that the $CMC$ level $S^{T_{cmc}}_{a}$ is convex. On the one hand, by \cite[Corollary~6.6]{Tamburelli}, the surface $S^{T_{cmc}}_{a}$ is obtained  by pushing to the future  the $k-$surface  $S^{T_{k}}_{-2a^{2}+2a\sqrt{a^{2}-1}+1}$ along its orthogonal geodesics for some positive time $t(a)>0$. On the other hand, since $M$ is future complete, the principal curvatures of $S^{T_{k}}_{-2a^{2}+2a\sqrt{a^{2}-1}+1}$ are negative. By \cite[Proposition~9.10]{BBZ2}, this implies that the principal curvatures of $S^{T_{cmc}}_{a}$ are also negative, which means exactly that $S^{T_{cmc}}_{a}$ is convex.
\end{proof}

Throughout the rest of the paper, we will assume that the $CMC$ takes its values in $\mathbb{R}_{+}^{*}$.

\section{Moncrief lines versus grafting rays}

In constant curvature spacetimes, convex Cauchy surfaces have a geometry comparable to that of cosmological levels. More precisely, we have the following theorem:
\begin{theo}\cite{mehdi3}
Let $S$ be a $C^{2}$ convex surface in a $2+1$ future complete $MGHC$ space-time $\left(M,g\right) $ of non negative constant curvature. Let $g_{S}$ be the Riemannian metric defined on $S$ by the restriction of the Lorentzian metric $g$. Then $\left(S,g_{S}\right)$ is $K^{4}-$bi-Lipschitz to $\left(S^{T_{cos}}_{\operatorname{Sup}_{S}T_{cos}}, g_{S^{T_{cos}}_{\operatorname{Sup}_{S}T_{cos}}} \right) $ where,
\begin{enumerate}
\item $K=\frac{\operatorname{Sup}_{S}T_{cos}}{\operatorname{Inf}_{S}T_{cos}}$ in the flat case;
\item $K=\frac{\sinh(\operatorname{Sup}_{S}T_{cos})}{\sinh(\operatorname{Inf}_{S}T_{cos})}$ in the de Sitter case.
\end{enumerate}
\end{theo}
Applying this to the $CMC$ time gives us:
\begin{prop}
\label{proppluimp}
Let $M$ be a $2+1$ non elementary future complete $MGHC$ space-time of non negative constant curvature. Let $T_{cos}$, $T_{cmc}$ be respectively the cosmological time and the $CMC$ time of $M$. Then for every $a\in \mathbb{R}_{+}^{*}$, the $CMC$ level $S^{T_{cmc}}_{a}$ is $K^{4}-$bi-Lipschitz to $S^{T_{cos}}_{a}$ where, 
\begin{enumerate}
\item $K=\frac{1}{2}$ in the flat case;
\item $K=2\cosh(\frac{a}{2})$ in the de Sitter case.
\end{enumerate}
\end{prop}

Let $S$ be a closed surface of hyperbolic type. Let $\lambda$ be a measured geodesic lamination on $S$. Let $M$ be the unique flat (respectively de Sitter) non elementary future complete $MGHC$ space-time of dimension $2+1$ associated to $\lambda$.  
Then Proposition \ref{proppluimp} together with \cite[Proposition~10.1]{mehdi3} yields:
\begin{cor}
Let $a\mapsto \left[g^{T_{cmc}}_{a}\right]$, $a\mapsto \operatorname{gra}_{\frac{\lambda}{a}}(S)$  be respectively the Moncrief line associated to $M$ and the grafting line defined by the measured geodesic lamination $(\lambda,\mu)$. Then, $\forall a\ll 1$, $$d_{\operatorname{Teich}}(\left[g^{T_{cmc}}_{a}\right],\operatorname{gra}_{\frac{\lambda}{a}}(S))\leq 4\log 3,$$ where $d_{\operatorname{Teich}}$ is the Teichm\"uller metric of $\operatorname{Teich}(S)$.
\end{cor}

\section{Convergence in Thurston's boundary of the Teichm\"uller space}
In this section, we will prove Theorem \ref{theoprinc}.\\ 
In fact, using Proposition \ref{proppluimp}, Theorem \ref{theoprinc} is a direct consequence of the following more general theorem:
\begin{theo}\label{theo2}
Let $S$ be a closed hyperbolic surface, and let $\lambda$ be a measured geodesic lamination on $S$. Let $h\in \operatorname{Teich}(S)$ and let $(h_{a})_{a \in \mathbb{R}^{*}{+}}$ be a curve in $\operatorname{Teich}(S)$, such that there exists a constant $K > 0$ verifying $d_{\operatorname{Teich}}(h_{a},gr_{\frac{\lambda}{a}}(h)) < K$, for every $a$ in $\mathbb{R}^{*}_{+}$. 
\begin{itemize}    
\item If the measured geodesic lamination $\lambda$ is uniquely ergodic, then the curve $(h_{a})_{a \in \mathbb{R}^{*}_{+}}$ converges to the projective lamination $ \left [ \lambda  \right ]$ on the Thurston boundary of the Teichm\"uller space $Teich(S)$. 
\item If $\lambda = \sum c_{i} \gamma_{i}$ is simplicial, then the curve  $(h_{a})_{a \in \mathbb{R}^{*}{+}}$ converges to the projective lamination $\left[\sum \gamma_{i}\right]$ on  the Thurston boundary of the Teichm\"uller space $Teich(S)$.
\end{itemize}
\end{theo}
This is a generalization of the principal theorem of Diaz and Kim \cite{Kim} initially obtained for grafting. Let $S$ be a closed hyperbolic surface and let $h\in \operatorname{Teich}(S)$. We begin with the following proposition:

\begin{prop}
\label{prop5}
Let $\left\{\gamma_{1},...,\gamma_{k} \right\}$ be a system of disjoint closed curve on $S$, and let $\lambda = \sum_{i} c_{i}\gamma_{i}$ be a simplicial measured geodesic lamination on $S$. Let $(h_{a})_{a \in \mathbb{R}^{*}_{+}}$ be a curve in $\operatorname{Teich}(S)$,  such that there exists  $K > 0$ satisfying $d_{Teich}(h_a,gr_{\frac{\lambda}{a}}) < K $ for every $a \in \mathbb{R}^{*}_{+}$. Then for each $i$ in $ \left\{1,...,k \right\}$ we have, 
$$\frac{2 e^{-2K} \theta}{2\theta + ca}\ell_{h}(\gamma_{i}) \leq \ell_{h_{a}}(\gamma_{i}) \leq \frac{e^{2K}\pi}{ \pi + c_{i}a}\ell_{h}(\gamma_{i}),$$
where $c = \max(c_{1},...,c_{k})$ and $\theta \in (0,\frac{\pi}{2})$ is a constant depending on the tuple of lengths $(\ell_{h}(\gamma_{1}),...,\ell_{h}(\gamma_{k}))$
\end{prop}
\begin{proof}
This follows directly from Proposition \ref{proppropprop} and \cite[Proposition~3.4]{Kim}
\end{proof}

Given a pants decomposition \(\{\gamma_{1}, \ldots, \gamma_{N}\}\), a curve \(\delta_{i}\) is said to be dual to \(\gamma_{i}\) if the intersection number \(i(\delta_{i}, \gamma_{j})\) is 0 for any \(j \neq i\), and \(i(\delta_{i}, \gamma_{i})\) is $1$ or $2$, depending on whether the two pants pieces of the decomposition glued along \(\gamma_{i}\) are the same or not, respectively. Before proceeding to the proof of Theorem \ref{theo2}, we need the following proposition:

\begin{prop}
\label{prop6}
 Let $\left\{\gamma_{1},...,\gamma_{k} \right\}$ be a set of simple closed curves and let $\lambda = \sum c_1\gamma_1 + ... + c_k\gamma_k$ be a simplicial measured geodesic lamination on $S$.  Let $(h_{a})_{a \in \mathbb{R}^{*}_{+}}$ be a curve in $\operatorname{Teich}(S)$,  such that there exists  $K > 0$ satisfying $d_{Teich}(h_a,gr_{\frac{\lambda}{a}}) < K $ for every $a \in \mathbb{R}^{*}_{+}$. We add curves $\gamma_{k+1}, \ldots, \gamma_{N}$ to form a system of pants decomposition $\left\{\gamma_{1},...,\gamma_{N} \right\}$. Let $\delta_{i}$ be the dual curve of $\gamma_{i}$. Then for every $i$ in $ \left\{1,...,N \right\}$, Then $\ell_{h_a}(\gamma_{i})Tw_{h_a}(\delta_{i},\gamma_{i})$ is bounded when $a\to +\infty$.
\end{prop}
\begin{proof}
This also follows immedialtely from \cite[Proposition~3.5]{Kim} and Lemma \ref{bounded}.
\end{proof}

\textbf{Proof of Theorem  \ref{theo2}}.
\textbf{$1.)$} Assume first that $\lambda$ is a maximal uniquely ergodic geodesic lamination. Let $(h_{a_{n}})_{n}$ be a subsequence of $(h_{a})_{a}$ converging to $[\mu]$. By contradiction, assume that the support of $\lambda$ and  $\mu$ are different. On the one hand, by maximality of $\lambda$ we must have $i(\lambda,\mu)\neq 0$. This implies that $\ell_{h_{a_{n}}}(\lambda)$ is not bounded. On the other hand, by Proposition \ref{proppropprop} and \cite[Corollary~3.2]{Kim}, we have $\ell_{h_{a_{n}}}(\lambda)\leq e^{2K}\ell_{gr_{\frac{\lambda}{a_{n}}}(h)}(\lambda)\leq e^{2K}\ell_{h}(\lambda) $ which is a contradiction. Now since $\lambda$ is uniquely ergodic, we have $[\lambda]=[\mu]$.\\
\textbf{$2.)$} Assume that $\lambda=\sum c_{i}\gamma_{i}$ is simplicial. Extend $\left\{\gamma_{1},...,\gamma_{k} \right\}$ to a pants decomposition $\left\{\gamma_{1},...,\gamma_{N} \right\}$. By \cite[Corollary~3.2]{Kim} and Proposition \ref{proppropprop}, the lengths $\ell_{h_a}(\gamma_{i})$ are all bounded independently of $a$. Then, \cite[Proposition 2.2]{Kim} gives us that for every simple closed curve $\beta$:
$$\ell_{h_a}(\beta) = \sum_{j = 1}^{N} i(\beta,\gamma_{j})\left( 2\ln\frac{1}{\ell_{h_a}(\gamma_{j})} + Tw_{h_a}(\beta,\gamma_{j})\ell_{h_a}(\gamma_{j})\right)  + O(1),$$
where $O(1)$ is a bounded function on $\beta$.

Now, on the one hand, by Proposition \ref{prop6}, all the terms  $Tw_{h_t}(\beta,\gamma_{i})\ell_{h_t}(\gamma_{i})$ are bounded  independetly of $a$. Moreover, all the terms $\ln\frac{1}{\ell_{h_a}(\gamma_{i})}$, for $i\in \left\{k+1,...,N \right\}$, are bounded from above. On the other hand, by Proposition \ref{prop5} we have that for every $i\in \left\{1,...,k \right\}$, $\frac{\ln\frac{1}{\ell_{h_a}(\gamma_{i})}}{\ln(a)}\longrightarrow 1$. Putting all this together give us $$\frac{1}{2\ln(a)}\ell_{h_a}(\beta)\longrightarrow \sum i(\beta,\gamma_{i})$$

\nocite{*}

\bibliographystyle{alpha}
\bibliography{bib.bib}
\end{document}